\documentclass[11pt]{amsart}
\usepackage{amssymb, latexsym, amsmath, amsfonts, tikz, amscd}
\usepackage{tikz-cd} 
\usepackage{hyperref, enumitem, fancyhdr}
\usepackage{color}
\newtheorem{thm}{Theorem}[section]
\newtheorem{cor}[thm]{Corollary}

\theoremstyle{definition}

\theoremstyle{remark}

\numberwithin{equation}{section}
\theoremstyle{remark}

\newcommand{\no}{\noindent}

\begin{document}
\title{On the group of automorphisms of the escaping sets of H\'{e}non maps}
\author{Mahima and Ratna Pal}

\address{M: Indian Institute of Science Education and Research Mohali, Knowledge City, Sector -81, Mohali, Punjab-140306, India}
\email{ph22080@iisermohali.ac.in, mahimamath06@gmail.com}
\address{RP: Theoretical Statistics and Mathematics Unit,
Indian Statistical Institute, Bangalore Centre, 8th Mile Mysore Road, Bangalore, 560059.}
\email{ratna.math@gmail.com, ratnapal@isibang.ac.in}

\begin{abstract}
In this article, we provide an explicit description of $\operatorname{Aut}_1(I_H^+)$, the group of automorphisms of the escaping set $I_H^+$ of a Hénon map that induce the identity on its fundamental group. More precisely, for a Hénon map
$
H(x,y)=(y,p(y)-\delta x),
$
where $p$ is a polynomial of degree $d\geq 2$ and $\delta\neq 0$, we show that $\operatorname{Aut}_1(I_H^+)$ is completely determined by the degree and the coefficients of $p$. Thus, the automorphism group of the escaping set, subject to the prescribed topological constraint, can be explicitly characterized in terms of the the underlying Hénon map.
	
\end{abstract}
\maketitle

\section{Introduction}
The group of automorphisms of a domain  $\Omega\subseteq \mathbb{C}^k$, for $k\geq 1$, is the group ${\rm{Aut}}(\Omega)$ of biholomorphic self maps of $\Omega$. The study of automorphism groups of complex domains is a central theme in complex analysis and several complex variables. Since biholomorphic domains have naturally isomorphic automorphism groups, the fundamental question is to understand to what extent ${\rm{Aut}}(\Omega)$ determines or reflects the complex geometry of the domains.  In particular, one would like to know which geometric, topological, or dynamical features of a domain are encoded in its group of biholomorphic self-maps. For bounded domains, the theory is particularly rich: the automorphism group is a real Lie group, and the interaction between its global structure, invariant metrics, and the boundary geometry of the domain has led to a number of fundamental rigidity results (\cite{Wong}, \cite{Rosay}, \cite{Kim}, \cite{Bedford and Dadok}, \cite{Bedford and Pinchuk}, \cite{Fridman and Poletsky}, \cite{Verma and Seshadri} and the references therein). In contrast, the study of automorphism groups of unbounded domains in higher dimensions is comparatively less explored. Unlike bounded domains, unbounded domains possess an additional geometric feature—the point at infinity—which can significantly influence both their intrinsic geometry and their automorphism groups. In particular, sequences of automorphisms may escape to infinity, making the analysis of orbit accumulation and boundary geometry more subtle.

\medskip 
\noindent
In the present article, we investigate the automorphism groups of certain dynamically defined unbounded domains associated with H\'{e}non maps in $\mathbb{C}^2$. Our primary objective is to understand to what extent the automorphism groups of these domains are determined by the H\'{e}non maps and reflect the underlying dynamics of the H\'{e}non map.  H\'{e}non maps are the polynomial automorphisms of $\mathbb{C}^2$ of the form 
\begin{equation}\label{form Henon}
H(x,y)=(y,p(y)-\delta x),
\end{equation} 
where $p$ is a polynomial of degree $d\geq 2$ and $\delta\neq 0$. These maps constitute the most important class of polynomial automorphisms of $\mathbb{C}^2$ exhibiting genuinely non-trivial dynamics. A fundamental result in the classification of polynomial automorphisms of $\mathbb{C}^2$ asserts that every polynomial automorphism of dynamical degree greater than or equal to $2$ is, up to polynomial conjugacy, a finite composition of H\'{e}non maps (\cite{FM}, \cite{Cantat and Dujardin}). Consequently, H\'{e}non maps provide a fundamental and sufficiently rich class of examples for studying the interaction between holomorphic dynamics and the the complex geometry of unbounded domains.

\medskip 
\noindent
The Escaping set of a H\'{e}non map $H$ is defined as
$$
I_H^+
=
\left\{z\in\mathbb C^2:\|H^n(z)\|\rightarrow\infty
\text{ as }n\rightarrow\infty\right\}.
$$
It turns out that $I_H^+$ is an $H$-invariant open connected pseudo-convex set.  An alternative and useful description of escaping sets comes from the Green's function $G_H^+$ associated to the H\'{e}non map $H$. The Green's function associated to $H$ is defined as
$$
G_H^+(z)
=
\lim_{n\to\infty}\frac{1}{d^n}\log^+\|H^n(z)\|,
$$
where $\log^+(x)=\max \{\log x,0\}$.  It turns out that 
\[
I_H^+=\{(x,y)\in \mathbb{C}^2: G_H^+(x,y)>0\}.
\]
Further, $G_H^+$ is pluri-harmonic on $I_H^+$.  Thus the Green's function serves as a natural tool for describing the analytic structure of $I_H^+$.  The complex analytic structure of $I_H^+$ was studied in detail by Hubbard and Oberste-Vorth (\cite{HOV}, \cite{Favre}). In particular, they established the fundamental group $\pi_1(I_H^+)\cong \mathbb{Z}[1/d].$
They further showed that $I_H^+$ admits a natural quotient representation
$$
I_H^+\cong {\mathbb{C}\times\left(\mathbb{C}\setminus\overline{\mathbb{D}}\right)}
/{\Gamma},
$$
 where $\Gamma$ is a discrete subgroup of
 $
\rm{Aut}\left(
\mathbb{C}\times
\left(\mathbb{C}\setminus\overline{\mathbb{D}}\right)
\right)
$
isomorphic to $\mathbb{Z}[1/d]/\mathbb{Z}$. This description provides an explicit model for the complex structure of the escaping set and, in particular, provides a natural framework for studying the group $\mathrm{Aut}_1(I_H^+)$,  the group of automorphisms of $I_H^+$ which lifts as an automorphism to the covering $\mathbb{C}\times\mathbb{C}\backslash\mathbb{\overline{D}}$. We prove the following:

\begin{thm}\label{main thm}
    Let $H(x,y)=(y,p(y)-\delta x)$ be a H\'{e}non map of degree $d$ as in (\ref{form Henon}) and let $p(y)=y^d+a_{d-2} y^{d-2}+\cdots+a_1 y+a_0$. Let $I_H^+$ be the escaping set of $H$. Let $i_p=gcd\{d^2-1, i:a_{d-i}\neq 0 \}$. Then $\mathrm{Aut}_1(I_H^+)\cong \mathbb{Z}_{i_p}\ltimes \mathbb{C}$.
 \end{thm}
 In a similar spirit as in Theorem \ref{main thm}, we can also give an explicit description of biholomorphisms between escaping sets of a pair of H\'{e}non maps $H$ and $F$ that can be lifted as an automorphism of $\mathbb{C}\times (\mathbb{C}\setminus \overline{D})$. We denote this set by $\mathrm{Biholo}_1(I_H^+\to I_F^+)$. Let $f_p:\mathbb{C}^2\to \mathbb{C}^2$ defined by $(x,y) \mapsto (e^{2d\pi i/i_p} x, e^{2\pi i/i_p} y)$. It turns out that the map $f_p$ is an automorphism of $I_H^+$ which generates the cyclic component $\mathbb{Z}_{i_p}$ of $\mathrm{Aut}_1(I_H^+)\cong \mathbb{Z}_{i_p}\ltimes\mathbb{C}$ (see Section 3). 
 
 \begin{thm}\label{Biholo escaping}
  Let $H$ and $F$ be two H\'{e}non maps of same degree $d$ such that $\mathrm{Biholo}_1(I_H^+\to I_F^+)$ is non-empty. Let $i_p=gcd\{d^2-1,i: a^H_{d-i}\neq 0\}$. Fix $f_\circ \in\mathrm{Biholo}_1(I_H^+\to I_F^+)$. Then for every $f\in \mathrm{Biholo}_1(I_H^+\to I_F^+)$, there exist an $r\in\{1,2,\dots,i_p\}$ and an automorphism $f_\gamma\in \mathrm{Aut}_1(I_H^+)$ with a lift of the form $(s,t)\mapsto (s+\gamma,t)$ of the covering space  $\mathbb{C}\times\mathbb{C}\setminus\mathbb{\overline{D}}$, where $\gamma\in \mathbb{C}$,  such that 
\[
f=f_\circ \circ f_p^r\circ f_\gamma.
\]
\end{thm}
 \noindent
 Now we record the following direct corollaries of Theorem \ref{main thm}. We omit proofs of Corollary \ref{pyd} and Corollary \ref{a0neq0} since they are particular cases of Theorem \ref{main thm}. 
  \begin{cor}\label{pyd}
Let $H$ be as in (\ref{form Henon}). Then $ \mathrm{Aut}_1(I_H^+)\cong \mathbb{Z}_{d^2-1}\ltimes \mathbb{C}$ if and only if $p(y)=y^d$. 
 \end{cor}
 
 \begin{cor}\label{a0neq0}
   Let $H$ be as in (\ref{form Henon}). Then   $\mathrm{Aut}_1(I_H^+)\cong \mathbb{C}$ if $a_0\neq 0$.
 \end{cor}
 
 \begin{cor}\label{centralizer of H}
    Let $H$ be a H\'enon map as in \ref{main thm}. Then centralizer of $H$ in $\mathrm{Aut}_1(I_H^+)$
    \[
C_{\mathrm{Aut}_1(I_H^+)}(H)\cong
\begin{cases}
\mathbb{Z}_{gcd(d-1,i_p)}\ltimes\mathbb{C}
& d=\delta, \\[6pt]
\mathbb{Z}_{gcd(d-1,i_p)}
& d\neq\delta.
\end{cases}
\]
\end{cor}
 \noindent 
 It follows from Theorem 3.2 in \cite{BV}, that $\mathbb{C}\leq \mathrm{Aut}_1(I_H^+) \leq \mathbb{Z}_{d^2-1}\ltimes \mathbb{C}$. While this result provides an upper bound for the finite cyclic component of $\operatorname{Aut}_1(I_H^+)$, it does not determine this component explicitly. The main purpose of the present article is to identify the finite cyclic subgroup precisely. In particular, Theorem \ref{main thm} gives an explicit description of $\operatorname{Aut}_1(I_H^+)$, thereby refining the structural result of \cite{BV}. Theorem \ref{main thm} also generalizes a result of Bousch in \cite{Bousch}, giving an explicit description of automorphism groups of escaping sets of quadratic H\'{e}non maps. 
 
 \medskip 
 \noindent
  It turns out that in some special cases $\mathrm{Aut}_1(I_H^+)$ gives complete description of $\mathrm{Aut}(I_H^+)$, for example when the degree of the H\'{e}non map is power of some prime (see \cite[Theorem 3.4]{BV}).  Further, it follows from \cite{Pal} (see Section 5) whenever the escaping sets of a H\'{e}non maps $F(x,y)=(y,p_F(y)-\delta_F x)$ of degree $d$  is  biholomorphic to the escaping set of $H$, then $F=A_1\circ H \circ  A_2$, for some diagonal linear maps $A_1$ and $A_2$ on $\mathbb{C}^2$. Thus we can conclude that escaping sets of a pair of H\'{e}non maps could possess the same automorphism groups inspite of not being biholomorphic to each other. In fact, we can further say:
\begin{cor}\label{continuum}
 There is a continuum of H\'{e}non maps with mutually non-biholomorphic escaping sets each having the same automorphism group. 
\end{cor}

\begin{cor}\label{global aut group} 
 Let $H$ be a H\'enon map as in \ref{main thm}. Then $\mathrm{Aut}_1(I_H^+)\cap \mathrm{Aut}(\mathbb{C}^2)\cong \mathbb{Z}_{i_p}$.
 \end{cor}
 
\medskip
\noindent
In one complex dimension, a simply connected domain with identically vanishing Kobayashi metric is necessarily biholomorphic to $\mathbb{C}$, by the Riemann mapping theorem. In higher dimensions, however, this phenomenon no longer true. A natural and important class of examples is given by the sublevel sets of the Green function associated with a H\'{e}non map $H$:
\[
\Omega_c=\{(x,y)\in \mathbb{C}^2: G_H^+(x,y)<c\}.
\]
Forn\ae ss showed in \cite{Fornaess} that these domains can be expressed as increasing unions of biholomorphic images of the unit ball, have identically vanishing Kobayashi metric, and are nevertheless not biholomorphic to $\mathbb{C}^2$. Such domains provide examples of so-called Short $\mathbb{C}^2$'s. For $c>0$, let  $\Omega_c'=\{(x,y)\in \mathbb{C}^2: 0<G_H^+<c\}$ be punctured Short $\mathbb{C}^2$.  Theorem \ref{main thm} suggests the following (also see \cite[Corollary 3.5]{BV}).
\begin{cor}\label{level sets of green}
 Let $H$ be a H\'enon map as in \ref{main thm}. Then $\mathrm{Aut}_1(\Omega'_c)\cong \mathbb{Z}_{i_p}\ltimes\mathbb{C}$, and $\mathrm{Aut}(\Omega_c)\cong \mathbb{Z}_{i_p}$.
\end{cor}

\begin{cor}\label{whole aut escaping}
     Let $H$ be a H\'enon map as in \ref{main thm} of degree $d=p_1^{m_1}p_2^{m_2}\dots p_l^{m_l}$ where for $1\leq i\leq l$, $p_i$ is prime. Then $\mathbb{Z}_{i_p}\ltimes\mathbb{C}\times\mathbb{Z}\leq \mathrm{Aut}(I_H^+)\leq \mathbb{Z}_{i_p}\ltimes\mathbb{C}\times\mathbb{Z}^l$.
\end{cor}

\begin{cor}\label{global whole aut group}
     Let $H$ be a H\'enon map as in \ref{main thm} of degree $d=p_1p_2\dots p_l$ where for $1\leq i\leq l$, $p_i$ is prime. Then $\mathrm{Aut}(I_H^+)\cap \mathrm{Aut}(\mathbb{C}^2)\cong \mathbb{Z}_{i_p}\times\mathbb{Z}$.
\end{cor}

\medskip
\no
We conclude this section by observing that, in contrast to the bounded setting, the automorphism groups of escaping sets and sub-level sets of Green functions may encode relatively limited information about the underlying dynamics and the intrinsic geometry of these domains. For instance, there exist quadratic H\'{e}non maps whose escaping sets have isomorphic automorphism groups, while exhibiting distinct geometric and dynamical features. In particular, Wolf \cite{Wolf} constructed quadratic H\'{e}non maps for which the boundaries of the corresponding escaping sets have different Hausdorff dimensions and these quadratic H\'{e}non maps can be chosen such that they have the same automorphism groups. Similarly, one can give examples of  quadratic H\'{e}non maps with the same automorphism group but different structures of their critical locus (see \cite{Lyubich and Robertson}). These examples illustrate that the automorphism group does not in general determine  finer geometric  dynamical properties of the underlying domains.

\medskip 
\noindent
{\textbf{Acknowledgement:}}  
The first author is supported by University Grants Commission (UGC). The second author is partially supported by the Mathematical Research Impact Centric Support\\
 (MTR/2023/001258) and the Advanced Research Grant (ARG/2025/007852) of the Anusandhan National Research Foundation (ANRF), India. A substantial part of this work was carried out during the first author's visit to the Indian Statistical Institute (ISI), Bangalore Centre. The first author sincerely thanks the Indian Statistical Institute, Bangalore Centre, for its warm hospitality and excellent research environment during her visit.

\section{Preliminaries}
Consider the H\'enon map $H$ of degree $d\geq 2$ defined in \eqref{form Henon} where $p(y)=y^d+a_{d-2}y^{d-2}+\cdots+a_1y+a_0$. The map $H$ extends as a birational map $\hat H:\mathbb{P}^2\to \mathbb{P}^2$ defined by
\[
\hat{H}[x:y:u]=[yu^{d-1}:y^d+a_{d-2}y^{d-2}u^2+\cdots+a_1yu^{d-1}+a_0u^d-\delta x u^{d-1}:u^d],
\] 
with lone indeterminacy point $[1:0:0].$ It turns out that $[0:1:0]$ is an attracting fixed point of $H$. This behavior at infinity provides the dynamical setting for the construction of the B\"ottcher coefficient.
For $R>0$, define
\begin{align*}
    V_R^+ & =\{(x,y)\in \mathbb{C}^2 : \max\{|x|,R\}<|y|\},\\
    V_R^- & =\{(x,y) \in \mathbb{C}^2 : \max\{|y|,R\}<|x|\},\\
    V_R & = \{(x,y)\in \mathbb{C}^2 : \max\{|x|, |y|\}\leq R\}.
\end{align*}
There exists sufficiently large $R>0$ such that $\{H^{-n}(V_R^+)\}_{n\geq 0}$ is a sequence of increasing open sets. Moreover, 
\[
I_H^+ =\bigcup_{n\geq 0} H^{- n}\left(V_R^+\right) \text{ and } K_H^+\subset V_R\cup V_R^-.
\]
These properties allow us to pass between the escaping set and the region $V_R^+$ where the B\"ottcher coordinate is defined.
For $n\geq 1$ and $(x,y)\in \mathbb{C}^2$, let $(x_n,y_n)=H^n(x,y).$ For $(x,y)\in V_R^+$, the B\"ottcher coordinate of $H$ is defined by
\[
\phi_H^+(x,y)=y\prod_{n=0}^{\infty}\left(\frac{y_{n+1}}{y_n^d}\right)^{1/d^{n+1}}.
\]
For sufficiently large $R$, this gives a holomorphic map $\phi_H^+:V_R^+\to \mathbb{C}\setminus\mathbb{\overline{D}}$ satisfying the functional equation 
\[
\phi_H^+\circ H=(\phi_H^+)^d.
\] 
Moreover,
\[
\phi_H^+(x,y)\sim y \text{ as $||(x,y)||\to \infty$ in } V_R^+,
\]
(see \cite{BS1,HOV,MNTU} for details).

\medskip
\noindent
The B\"ottcher coordinate does not extend as a single valued holomorphic function to all of $I_H^+$. It admits analytic continuation along every path in $I_H^+$ starting in $V_R^+$. The resulting multi-valued function can be made single valued on a natural covering space. The fundamental group of $I_H^+$ turns out to be $\mathbb{Z}[1/d]$. Further the fundamental group of $V_R^+$ is the subgroup $\mathbb{Z}$. We therefore consider the covering corresponding to this subgroup. We briefly recall the construction for the covering $\hat{I}_H^+$ of $I_H^+$ corresponding to the subgroup $\mathbb{Z}$, which turns out to be biholomorphic to $\mathbb{C}\times\mathbb{C}\setminus\mathbb{\overline{D}}$ (see \cite{HOV,MNTU}). Fix $z_0\in V_R^+$. Consider the set of pairs $(z,C)$, $z\in I_H^+$ where $C\subset I_H^+$ is a path from $z_0$ to $z$. Define an equivalence relation by
\[
(z,C)\sim(z',C') \iff z=z' \text{ and } [C\overline{C'}]\in \mathbb{Z}.
\]
Then $\hat{I}_H^+$ is the set of equivalence classes. The projection $\pi_H:\hat{I}_H^+\to I_H^+$ defined by $[(z,C)]\mapsto z$ is a covering map.
Let $\hat{V}_R^+=\{[(z,C)]\in\hat{I}_H^+ : z\in V_R^+\}$. Then $\pi_H|_{\hat{V}_R^+}:\hat{V}_R^+\to V_R^+$ is biholomorphic. On the set $\hat{V}_R^+$, the lifted B\"ottcher coordinate  agrees with the original function $\phi_H^+$.

\medskip 
\noindent
Now we briefly recall the construction of the additional coordinate which together with the B\"ottcher coordinate gives the biholomorphism from $\hat{I}_H^+$ to $\mathbb{C}\times\mathbb{C}\setminus\overline{\mathbb{D}}$ from \cite{MNTU}. For large constant $M>0$, let 
\[
U_R^+=\left\{(x,y)\in V_R^+ : \left\lvert \phi_H^+(x,y)\right\rvert>M\max\{R,|x|\}\right\}.
\]
It follows that $H(U_R^+)\subset U_R^+$.
Since $\phi_H^+(x,y)\sim y \text{ as $\lVert (x,y)\rVert \to \infty$ in } V_R^+$, the map 
\[
(x,y)\mapsto (x,t)=\left (x,\phi_H^+(x,y)\right)
\]
is a biholomorphism from $U_R^+$ onto its image. We write its inverse as
\[
(x,t)\mapsto \left(x,\lambda_H^+(x,t)\right),
\] 
In particular, $\phi_H^+\bigl(0,\lambda_H^+(0,t)\bigr)=t$. The asymptotic relation $\phi_H^+(0,y)\sim y$ gives a Laurent
expansion at infinity of the form
\[
\phi_H^+(0,y)=y+\frac{L_1^H}{y}+\frac{L_2^H}{y^2}+\frac{L_3^H}{y^3}+\cdots,
\]
where $L_j^H\in\mathbb{C}$. Its inverse therefore has an expansion
\[
\lambda_H^+(0,t)=t+D_0^H+\frac{D_1^H}{t}+\frac{D_2^H}{t^2}+\frac{D_3^H}{t^3}+\cdots.
\]
Substituting this expansion into $\phi_H^+\bigl(0,\lambda_H^+(0,t)\bigr)=t$ shows that $D_0^H=0$.
Consequently,
\[
\lambda_H^+(0,t)=t+\frac{D_1^H}{t}+\frac{D_2^H}{t^2}+\frac{D_3^H}{t^3}+\cdots.
\]

\medskip
\no
For $(x,y)\in U_R^+$, let $t=\phi_H^+(x,y)$. Define
\[
s_1(x,y)=\phi_H^+(x,y) \int_0^x \frac{\partial\lambda_H^+}{\partial t}(\xi,t)\mathrm{d}\xi ,
\]
In the $(s_1,t)$ coordinates, the H\'{e}non map $H$ takes the form
\[
(s_1,t)\mapsto \left(\frac{\delta}{d}s_1+\tilde{Q}_H(t),t^d\right),
\]
where $\tilde{Q}_H(t)$ is holomorphic for $|t|>RM$. More precisely, 
\[
\tilde{Q}_H(t)=t^d\int_0^{\lambda_H^+(0,t)}\frac{\partial\lambda_H^+}{\partial t}(\xi,t^d)\mathrm{d}\xi.
\]
Since $\lambda_H^+(0,t)-t$ remains bounded as $|t|\to \infty$, we have
$\left\lvert \tilde{Q}_H(t)\right \rvert\leq {\lvert t \rvert}^{d+1}$ for $|t|>RM$. Consequently, $\tilde{Q}_H$ has a Laurent expansion
\[
\tilde{Q}_H(t)=t^{d+1}+A^H_{d}t^d+A^H_{d-1}t^{d-1}+\cdots+A^H_0+\frac{c^H_1}{t}+\frac{c^H_2}{t^2}+\cdots,
\]
where $A^H_j,c^H_i\in \mathbb{C}$ for $0\leq j\leq d$ and $i\geq 1$. Further, the identity $\phi_H^+\left(0,\lambda^+_H(0,t)\right)=t$ forces $A^H_d=0$.
Let 
\[Q_H(t)=t^{d+1}+A^H_{d-1}t^{d-1}+A^H_{d-2}t^{d-2}+\cdots+A^H_0,\]
be the polynomial part of the Laurent expansion of $\tilde{Q}_H$. Note that $Q_H$ is polynomial part of $t^d\lambda_H^+(0,t)$. Define
\begin{align*}
R_H(t) & = \left(\frac{d}{\delta}\right)\left(\tilde{Q}_H(t)-Q_H(t)\right)+\left(\frac{d}{\delta}\right)^2\left(\tilde{Q}_H(t^d)-Q_H(t^d)\right)+\cdots.
\end{align*}
The series converges normally for sufficiently large $|t|$. 
Now define
\begin{equation}\label{psi_H^+}
s=\psi_H^+(x,y)=\phi_H^+(x,y) \int_0^x \frac{\partial\lambda_H^+}{\partial t}(\xi,t)\mathrm{d}\xi +R_H(t), \quad t=\phi_H^+(x,y).
\end{equation}
Then the map $H$ in the $(s,t)$ coordinates is given by
\[
(s,t)\mapsto \left(\frac{\delta}{d}s+Q_H(t),t^d\right).
\]
Since, for sufficiently large $R$, the map $(x,y)\mapsto(s_1,t)$ is a biholomorphism from $U_R^+$ onto its image, the map \[\Phi_H^+:U_R^+\to \mathbb{C}\times\mathbb{C}\backslash\mathbb{\overline{D}}\] defined by \[\Phi_H^+(x,y)=\left(\psi_H^+(x,y),\phi_H^+(x,y)\right),\] is injective holomorphic map on $U_R^+$. Furthermore,
\begin{equation}\label{psi equation}
    \psi_H^+\left(H(x,y)\right)=\left(\frac{\delta}{d}\right)\psi_H^+(x,y)+Q_H\left(\phi_H^+(x,y)\right).
\end{equation}

\medskip
\noindent
Let $\hat{H}:\hat I_H^+\to \hat I_H^+$ denote the lift of $H$. Let $\hat{U}_R^+=\{[(z,C)]\in \hat{V}_R^+ :  z \in U_R^+\}$. Then by \cite[Lemma 7.3.8]{MNTU}, 
\begin{equation}\label{hatU_H}
    \hat{U}_H^+=\bigcup_{n\geq 0} \hat{H}^{-n}(\hat{U}_R^+).
\end{equation}
For $[(z,C)]\in \hat{U}_R^+$, define \[\hat{\Phi}_H^+\left([(z,C)]\right)= \left(\psi_H^+(z),\phi_H^+(z)\right).\]
Although the functions on the right are defined only on $U_R^+$, the functional equations \eqref{hatU_H} and \eqref{psi equation} allows both coordinates to be continued consistently across appearing in \eqref{hatU_H}. Hence $\hat{\phi}_H^+$ extends to a holomorphic map 
\[\hat{\Phi}_H^+:\hat{I}_H^+\to \mathbb{C}\times\mathbb{C}\backslash\mathbb{\overline{D}}.\]
The resulting map is a biholomorphism (see \cite[Theorem 7.3.10]{MNTU}).
This analytic structure, introduced by Hubbard and Oberste-Vorth (\cite{HOV}), will be used below to study the group of automorphisms of $I_H^+$.

\section{Automorphism group of the escaping set}
In this section we provide proof of Theorem \ref{main thm}. The proof is divided into several steps. In the first step we prove $\mathrm{Aut}_1(I_H^+)\leq \mathbb{Z}_{i_p}\ltimes\mathbb{C}$, and in rest of the steps we prove $\mathbb{Z}_{i_p}\ltimes\mathbb{C}\leq \mathrm{Aut}_1(I_H^+)$. First we define an operation on $\mathbb{Z}_{i_p}\ltimes\mathbb{C}$. Consider the homomorphism $\theta:\mathbb{Z}_{i_p}\to \mathrm{Aut}(\mathbb{C})$ defined by $\alpha\mapsto \alpha^{\mathbb{C}}$, where $\alpha^{\mathbb{C}}\in \mathrm{Aut}(\mathbb{C})$ is defined by $\alpha^{\mathbb{C}}(x)=\alpha^{d+1}x$.
For $(\alpha_1,\gamma_1),(\alpha_2,\gamma_2)\in \mathbb{Z}_{i_p}\ltimes\mathbb{C}$ define 
\[
(\alpha_1,\gamma_1)*(\alpha_2,\gamma_2)=\left(\alpha_1\alpha_2,\gamma_1\alpha_1^{\mathbb{C}}(\gamma_2)\right).\] 

\medskip
\no
{\it Step 1.} 
It follows from \cite{Pal} and \cite{MP} that every $f \in \mathrm{Aut}_1(I_H^+)$ admits a lift $\tilde{f}\in\mathrm{Aut}_1(\mathbb{C}\times\mathbb{C}\setminus\overline{\mathbb{D}})$ of the form
\[
\tilde{f}(s,t) = (\beta_f s + \gamma_f,\ \alpha_f t),
\]
where $\beta_f = \alpha_f^{\,d+1},$ $\alpha_f^{\,d^2-1} = 1,$ $\gamma_f \in \mathbb{C}.$ The compatibility of this lift with the Böttcher coordinate of $H$ gives, as in \cite[Equation (5.17)]{Pal}, $\beta_f\, p(y) = \alpha_f\, p(\alpha_f y),$ $ \forall\, y \in \mathbb{C}$. Comparing the coefficient of $y^{d-i}$ on both sides gives
\[
\beta_f\, a_{d-i} = \alpha_f^{\,d-i+1}\, a_{d-i}, \qquad 2 \le i \le d.
\]
Substituting $\beta_f = \alpha_f^{\,d+1}$ yields $\alpha_f^{\,d-i+1}\big(\alpha_f^{\,i} - 1\big)\, a_{d-i} = 0$. Since $\alpha_f \neq 0$,
\begin{equation}\label{alpha coeff poly}
    a_{d-i} = \alpha_f^{\,i}\, a_{d-i}
\end{equation}
for $2\leq i\leq d$. Let 
\[
S=\{i\in\{2,3,\cdots,d\}:a_{d-i}\neq 0\}=\{i_1,i_2,\cdots ,i_{r_0}\}.
\]
It follows from \eqref{alpha coeff poly}, $a_{d-i}=\alpha_f^i a_{d-i}$ for $2\leq i\leq d$. Thus $\alpha_f^i=1$ for $i\in S$. Since $i_p=gcd\{d^2-1,i_1,i_2,\cdots,i_r\}$, there exist integers $n_1,n_2,\cdots,n_r$ such that $i_p=n_1i_1+n_2i_2+\cdots+n_ri_r$. Therefore 
\[
\alpha_f^{i_p}=1.
\] 
Following \cite[Section~3.1]{BV}, consider the map
\[
\mathcal{F}:\mathrm{Aut}_1(I_H^+)\to \mathbb{Z}_{i_p}\ltimes \mathbb{C}\quad \text{defined by }f\mapsto (\alpha_f,\gamma_f).
\]
We show that $\mathcal{F}$ is an injectice homomorphism.
Since other lifts of $f$ are $\tilde{f}\circ\gamma_{\frac{k}{d^n}}$ for $k\geq 1$, $n\geq0$, $\tilde{f}$ is unique lift of its form. Therefore, $\mathcal{F}$ is well defined.
Let $f_1,f_2\in \mathrm{Aut}_1(I_H^+)$ with $\tilde{f}_1(s,t)=(\beta_{f_1}s+\gamma_{f_1},\alpha_{f_1})$ and $\tilde{f}_2(s,t)=(\beta_{f_2}s+\gamma_{f_2},\alpha_{f_1})$. Then 
\[
\widetilde{f_1\circ f_2}(s,t)=\left(\beta_{f_1}\beta_{f_2}s+\beta_{f_1}\gamma_{f_2}+\gamma_{f_1},\alpha_{f_1}\alpha_{f_2}t\right).
\]
Thus 
\begin{align*}
\mathcal{F}(f_1\circ f_2)&=\left(\alpha_{f_1}\alpha_{f_2},\alpha^{d+1}_{f_1}\gamma_{f_2}+\gamma_{f_1}\right)\\
&=\mathcal{F}(f_1)*\mathcal{F}(f_2),
\end{align*}
which implies $\mathcal{F}$ is a group homomorphism. Now let $\mathcal{F}(f)=(1,0)$. Then lift of $f$, $\tilde{f}$, is given by $\tilde{f}(s,t)=(s,t)$. Thus $f$ is identity map. Therefore, 
\[\mathrm{Aut}_1(I_H^+)\cong \mathcal{F}(\mathbb{Z}_{i_p}\ltimes\mathbb{C})\leq \mathbb{Z}_{i_p}\ltimes\mathbb{C}.\]

\medskip
\no
{\it Step 2.} Let $j_p=e^{2\pi i/i_p}$. In this step, we prove that $f_p\in \mathrm{Aut}(\mathbb{C}^2)$ defined by $f_p(x,y)=(j_p^dx,j_py)$ is an automorphism of $I_H^+$. Let $g_p\in \mathrm{Aut}(\mathbb{C}^2)$ be defined by $g_p(x,y)=(j_px,j_p^dy)$. Consider
\[
p(j_py)=j_p^dy^d+a_{d-2}j_p^dj_p^{-2}y^{d-2}+a_{d-3}j_p^dj_p^{-3}y^{d-3}+\cdots+a_{1}j_p^dj_p^{-(d-1)}y+a_0j_p^dj_p^{-d}.
\]
Since $j_p^{i_p}=1$ and if $a_{d-i}\neq 0$ then $i_p$ divides $i$, whenever $a_{d-i}\neq 0$ then $j_p^i=1$. Thus,
\[
p(j_py)=j_p^dp(y).
\] 
Inductively, for $k\geq 1$, 
\begin{equation}\label{pjy}
    p(j_p^ky)=j_p^{d k}p(y).
\end{equation}
Note that $H\left(f_p(x,y)\right)=g_p\left(H(x,y)\right)$. Similarly, since $p(j_p^d y)=j_p p(y)$, $H\left(g_p(x,y)\right)=f_p\left(H(x,y)\right)$. Since $H^2\circ f_p=H\circ g_p\circ H=f_p\circ H^2$, $H^2$ and $f_p$ commutes. Inductively, for every $n\geq 1$, $H^{2n}\circ f_p=f_p\circ H^{2 n}$. Now
\begin{align*}
    z\in I_H^+ & \iff H^{2n}(z)\to \infty \text{ as } n\to\infty\\
    & \iff f_p \circ H^{2 n}(z)\to \infty \text{ as } n\to \infty\\
    & \iff H^{2 n}\left(f_p(z)\right)\to \infty \text{ as } n\to \infty\\
    & \iff f(z)\in I_H^+.
\end{align*}
Thus $f(I_H^+)=I_H^+$.

\medskip
\no
{\it Step 3.} In this step, we compute $\phi_H^+\left(f_p(x,y)\right)$ and $\psi_H^+\left(f_p(x,y)\right)$ for $(x,y)\in U_R^+$.
For $n\geq 1,$ let $(x_n,y_n)=H^n(x,y)$. First, using induction, we prove that $x_n\left(j_p^dx,j_py\right)=j_p^{d^{n-1}}x_n(x,y)$ and $y_n\left(j_p^dx,j_py\right)=j_p^{d^n}y_n(x,y)$. For $n=1$
\[
\left(x_1\left(j_p^dx,j_py\right),y_1\left(j_p^dx,j_py\right)\right)= H\left(j_p^dx,j_py\right)= \left(j_px_1(x,y),j_p^dy_1(x,y)\right).
\]
Let the result follows for $n-1$. Then for $n$,
\begin{align*}
     H^n(j_p^dx,j_py)&= H(H^{n-1}(j_p^dx,j_py))\\
     &=H(x_{n-1}(j_p^dx,j_py),y_{n-1}(j_p^dx,j_py))\\
     &=H(j_p^{d^{n-2}}x_{n-1}(x,y),j_p^{d^{n-1}}y_{n-1}(x,y))\\
     &=(j_p^{d^{n-1}}y_{n-1}(x,y),p(j_p^{d^{n-1}}y_{n-1}(x,y))+\delta j_p^{d^{n-2}}x_{n-1}(x,y))\\
     &=(j_p^{d^{n-1}}x_n(x,y),j_p^{d^n}p(y_{n-1}(x,y))+\delta j_p^{d^n}x_{n-1}(x,y))\\
     &=(j_p^{d^{n-1}}x_n(x,y),j_p^{d^n}y_n(x,y)).
\end{align*}
Thus, 
\[
\frac{y_{n}\left(j_p^dx,j_py\right)}{\left(y_{n-1}\left(j_p^dx,j_py\right)\right)^d}=\frac{j_p^{d^n}y_n(x,y)}{\left(j_p^{d^{n-1}}y_{n-1}(x,y)\right)^d}=\frac{y_n(x,y)}{y_n(x,y)^d},
\]
therefore,
\begin{equation}\label{phi f}
    \phi_H^+\left(j_p^dx,j_py\right)=j_py\cdot \left(\frac{y_1}{y^d}\right)^{\frac{1}{d}}. \cdots . \left(\frac{y_{n}}{y_{n-1}^d}\right)^{\frac{1}{d^{n+1}}}.\cdots=j_p\phi_H^+(x,y).
\end{equation}
Since $\phi_H^+\left(x,\lambda_H^+(x,t)\right)=t$ on $\left\{(x,t)\in \mathbb{C}^2:|t|>M\max\{R,|x|\}\right\} $, $\phi_H^+\left(j_p^d x,\lambda_H^+(j_p^d x,j_pt)\right)= j_pt$. Further, by \eqref{phi f}, it follows that $\phi_H^+\left(j_p^dx,j_p\lambda_H^+(x,t)\right)=j_p\phi_H^+\left(x,\lambda_H^+(x,t)\right)=j_pt$. Therefore,
$\phi_H^+\left(j_p^dx,j_p\lambda_H^+(x,t)\right)=\phi_H^+\left(j_p^dx,\lambda_H^+(j_p^dx,j_pt)\right)$. Since $\phi_H^+(j_p^dx,\cdot)$ is injective on $U_R^+$, $\lambda_H^+(j_p^dx,j_pt)=j_p\lambda_H^+(x,t)$, which further gives
\[
\frac{\partial\lambda_H^+}{\partial t}(j_p^dx,j_pt)=\frac{\partial\lambda_H^+}{\partial t}(x,t).
\]
Similarly, since $p(j_p^dy)=j_p^{d^2}p(y)=j_pp(y)$, we have for $n=1$, $x_1(j_px,j_p^dy)=j_p^dx_1(x,y)$ and $y_1(j_px,j_p^dy)=j_py_1(x,y)$. Inductively, for $n\geq 1$,
\[
x_n(j_px,j_p^dy)=j_p^{d^n}x_n(x,y), \quad y_n(j_px,j_p^dy)=j_p^{d^{n-1}}y_n(x,y).
\]
Thus, \[\phi_H^+\left(j_px,j_p^dy\right)=j_p^d\phi_H^+(x,y),\] which further implies
\[
\lambda_H^+(j_px,j_p^dy)=t^d \lambda_H^+(x,y) \text{ and } \frac{\partial\lambda_H^+}{\partial t}(j_px,j_p^dt)=\frac{\partial\lambda_H^+}{\partial t}(x,t).
\]
Now,
\begin{align*}
s_1\left(j_p^dx,j_pt\right)&=j_pt\int_0^{j_p^dx}\frac{\partial\lambda_H^+}{\partial t}(\xi,j_pt)\mathrm{d}\xi\\
& = j_p^{d+1}t\int_0^x\frac{\partial\lambda_H^+}{\partial t}(j_p^d\xi,j_pt)\mathrm{d}\xi\\
&= j_p^{d+1}t\int_0^x\frac{\partial\lambda_H^+}{\partial t}(\xi,t)\mathrm{d}\xi\\
&=j_p^{d+1}s_1(x,y).
\end{align*}
Further,
\begin{align*}
    \tilde{Q}_H(j_pt)& = j_p^dt^d\int_0^{\lambda_H^+(0,j_pt)}\frac{\partial \lambda_H^+}{\partial t}(\xi,j_p^dt^d)\mathrm{d}\xi \\
    &= j_p^dt^d\int_0^{j_p\lambda_H^+(0,t)}\frac{\partial \lambda_H^+}{\partial t}(\xi,j_p^dt^d)\mathrm{d}\xi \\
    &= j_p^dt^d\int_0^{\lambda_H^+(0,t)}\frac{\partial \lambda_H^+}{\partial t}(j_p\xi,j_p^dt^d)\mathrm{d}\xi\\
    &=j_p^{d+1}\tilde{Q}_H(t).
\end{align*}
Inductively, for $n\geq 1$, $\tilde{Q}_H\left(j_p^{d^n} t\right)=j_p^{(d+1)d^{n}}\tilde{Q}_H(t)=j_p^{d^{n+1}}j_p^{d^n}\tilde{Q}_H(t)$. Now, $j_p^{d^2}=j_p$, therefore, if $n$ is even then $j_p^{d^n}=j_p$ and if $n$ is odd then $j_p^{d^n}=j_p^d$. Thus $\tilde{Q}_H\left(j_p^{d^n} t\right)=j_p^{d+1}\tilde{Q}_H$. Since $Q_H$ is polynomial part of $\tilde{Q}_H$, for every $n\geq 1$, $Q_H\left(j_p^{d^n}t\right)=j_p^{d+1} Q_H(t)$. Now,
\[
R_H(j_pt)=\left(\frac{d}{\delta}\right)\left(\tilde Q_H(j_pt)-Q_H(j_pt)\right)+\left(\frac{d}{\delta}\right)^2\left(\tilde Q_H\left(j_p^dt^d\right)-Q_H\left(j_p^dt^d\right)\right)+\cdots= j_p^{d+1}R_H(t).
\]
Therefore 
\begin{equation}\label{psi f}
\begin{split}
    \psi_H^+\left(j_p^dx,j_p y\right) & = s_1\left(j_p^dx,\phi_H^+\left(j_p^dx,j_py\right)\right)+R_H\left(\phi_H^+\left(j_p^dx,j_py\right)\right)\\
    & = s_1\left(j_p^d x,j_p\phi_H^+(x,y)\right)+R_H\left(j_p\phi_H^+(x,y)\right)\\
    & = j_p^{d+1}s_1\left(x,\phi_H^+(x,y)\right)+j_p^{d+1}R_H\left(\phi_H^+(x,y)\right)\\
    & = j_p^{d+1}\psi_H^+(x,y).
\end{split}
\end{equation}

\medskip
\no
{\it Step 4.} In this step, we see that $\tilde{f}_p:\mathbb{C}\times\mathbb{C}\setminus\mathbb{\overline{D}}\to\mathbb{C}\times\mathbb{C}\setminus\mathbb{\overline{D}}$ defined by $\tilde{f}_p(s,t)=\left(j_p^{d+1}s,j_pt\right)$ is a lift of $f_p$. To prove this, we need to prove that $\pi_H\circ(\hat{\Phi}_H^+)^{-1}\circ \tilde{f}_p=f_p\circ \pi_H\circ (\hat{\Phi}_H^+)^{-1}$ on $\mathbb{C}\times\mathbb{C}\setminus\mathbb{\overline{D}}$. Since $\pi_H:\hat{U}_R^+\to U_R^+$ is a biholomorphism, it is enough to show that 
\begin{equation}\label{lift identity}
    \hat{\Phi}_H^+\circ\pi_H^{-1}\circ f_p = \tilde{f}_p\circ \hat{\Phi}_H^+\circ \pi_H^{-1}
\end{equation}
on $U_R^+$. Then the result follows from identity theorem. Let $W_R^+=\Phi_H^+(U_R^+)=\hat{\Phi}_H^+(\hat{U}_R^+)$. Observe that the map $\hat{\Phi}_H^+\circ\pi_H^{-1}:U_R^+\to W_R^+$ is defined by 
\[
(x,y)\mapsto\left(\psi_H^+(x,y),\phi_H^+(x,y)\right).
\]
It follows from \eqref{phi f} and \eqref{psi f} that
\[
\left(\psi_H^+\circ f_p(x,y),\phi_H^+\circ f_p(x,y)\right)=\tilde{f}_p\left(\psi_H^+(x,y),\phi_H^+(x,y)\right).
\]
Thus \eqref{lift identity} holds.

\medskip
\no
{\it Step 5.} In this step, similar to \cite[Section~3.1]{BV}, we prove that for every $\gamma\in\mathbb{C}$ there exists an $f\in \mathrm{Aut}_1(I_H^+)$ which lifts to the map $\tilde{f}(s,t)=(s+\gamma,t)$. We define $f$ using $\pi_H\circ \tilde{f}=f\circ\pi_H$. Let $z=(x,y)\in I_H^+$, choose $\tilde{z}\in \pi_H^{-1}(z)$. Define $f(z)=\pi_H\circ\tilde{f}(\tilde{z})$. We first verify that $f$ is well defined. Let $(s_1,t_1),(s_2,t_2)\in \mathbb{C}\times\mathbb{C}\setminus\mathbb{\overline{D}}$ be such that $\pi_H(s_1,t_1)=z=\pi_H(s_2,t_2)$. Then there exists $k/d^n\in \mathbb{Z}[1/d]/\mathbb{Z}$ such that $(s_1,t_1)=\gamma_{\frac{k}{d^n}}(s_2,t_2)$. It is enough to show that $\pi_H\circ\tilde{f}(s_1,t_1)=\pi_H\circ\tilde{f}(s_2,t_2)$, that is, $\pi_H\circ\tilde{f}\left(\gamma_{\frac{k}{d^n}}(s_2,t_2)\right)=\pi_H\circ\tilde{f}(s_2,t_2)$, that is,
     \[
      \pi_H\left(s_2+ \frac{d}{\delta} \sum_{l=0}^{n-1} {\left(\frac{d}{\delta}\right)}^l\left(Q_H(t_2^{d^l})-Q_H\left( {\left(e^{\frac{2 k\pi i}{d^n}}t_2)\right)}^{d^l} \right)\right)+\gamma,
 e^{\frac{2 k\pi i}{d^n}}t_2\right)=\pi_H(s_2+\gamma,t_2),
     \]
     that is,
     \[
     \pi_H\circ\gamma_{\frac{k}{d^n}}(s_2+\gamma,t_2)=\pi_H(s_2+\gamma,t_2),
     \]
     which is true. Hence $f$ is well defined. Finally, observe that every other lift of $f$ is of the form $\tilde{f}\circ \gamma_{k/d^n}$ for some $k\geq 0, n\geq 1$. For $k/d^n\neq 0 \mod \mathbb{Z}$, the first coordinate of $\tilde{f}\circ \gamma_{k/d^n}(s,t)$ is of the form $s+\gamma+\eta_{\frac{k}{d^n}}(t)$ where $\eta_{\frac{k}{d^n}}(t)$ is non constant. Consequently, no lift of $f$ other than $\tilde{f}$ itself can be of the form $(s,t)\mapsto(s+\gamma',t)$ for a constant $\gamma'\in \mathbb{C}$. Thus $\mathbb{C}\leq \mathrm{Aut}_1(I_H^+)$.

\medskip
\no
Since $\mathrm{Aut}_1(I_H^+)$ has an element of order $i_p$ and $\mathbb{C}\leq \mathrm{Aut}_1(I_H^+)\leq \mathbb{Z}_{i_p}\ltimes\mathbb{C}$, $\mathbb{Z}_{i_p}\ltimes\mathbb{C}\leq \mathrm{Aut}_1(I_H^+)$.

\subsection*{Proof of Corollary \ref{centralizer of H}}
For every $\gamma\in \mathbb{C}$, let $f_\gamma\in \mathrm{Aut}_1(I_H^+)$ be the automorphism whose lift is the map $\tilde{f}_\gamma(s,t)=(s+\gamma,t)$. Then $\tilde{H}\circ \tilde{f}_\gamma=\tilde{f}_{(\delta/d)\gamma}\circ\tilde{H}$, which implies 
\begin{equation}\label{f gamma commute}
    H\circ f_\gamma=f_{(\delta/d)\gamma}\circ H.
\end{equation}
Moreover, $H\circ f_p=g_p\circ H$. Snce $g_p=f_p^d$, we obtain 
\begin{equation}\label{Hcirc fp commute}
H\circ f_p=f_p^d\circ H.
\end{equation}
Let $f\in\mathrm{Aut}_1(I_H^+)$. Then $f$ admits a lift of the form $\tilde{f}(s,t)=(\alpha^{d+1}s+\gamma,\alpha t)$ where $\alpha^{i_p}=1$ and $\gamma\in \mathbb{C}$. Since $j_p=e^{2\pi i/i_p}$, there exists $r\in\{1,2,\dots,i_p\}$ such that $\alpha=j_p^r$. Consequently, $\tilde{f}=\tilde{f}_p^r\circ \tilde{f}_{\alpha^{-(d+1)}\gamma}$ and therefore
\begin{equation}\label{element in Aut escaping}
f=f_p^r\circ f_{\gamma'},
\end{equation}
for some $\gamma'\in \mathbb{c}$.
We first determine when $f_\gamma$ commute with $H$. By \eqref{f gamma commute}, for $\gamma\neq 0$, $f_\gamma$ commute with with $H$ if and only if $d=\delta$. Now, we determine the elements in the cyclic subgroup generated by $f_p$ that commute with $H$. Iterating \eqref{Hcirc fp commute} gives $H\circ f_p^k=f_p^{kd}\circ H$. Since $f_p$ has order $i_p$,
\[
H\circ f_p^k=f_p^k\circ H \iff f_p^k= f_p^{kd}\iff f_p^{(d-1)k}= Id \iff (d-1)k=0\mod i_p.
\]
Fo a general element $f=f_p^r\circ f_\gamma$,
\[H\circ f=H\circ f_p^r\circ f_{\gamma'}= f_p^{rd}\circ f_{(\delta/d)\gamma'}\circ H,\]
and
\[
f\circ H=  f_p^r\circ f_{\gamma'}\circ H.
\]
Thus $f$ commutes with $H$ if and only is $(d-1)r=0\mod i_p$ and $(\delta/d)\gamma'=\gamma'$.

\subsection*{Proof of Corollary \ref{continuum}} 
For each $c\in\mathbb{R}_{>0}$, consider the H\'enon map
\[
H_c(x,y)=(y,p_c(y)-\delta x),
\]
where $p_c(y)=y^d+c$. Since $c\neq 0$, the only nonzero
coefficient $a_{d-i}$ with $2\leq i\leq d$ occurs for
$i=d$. Hence $i_{p_c}=\gcd\{d^2-1,d\}=1$.
It follows from \ref{main thm} that $\mathrm{Aut}_1(I_{H_c}^{+})
\cong \mathbb{C}$, for every $c\in\mathbb{R}_{>0}$.
We now prove that the escaping sets $I_{H_c}^{+}$ are
pairwise non-biholomorphic. Suppose that, for some
$c,c'\in\mathbb{R}_{>0}$, there exists a biholomorphism from $I_{H_c}^{+}$ to $I_{H_{c'}}^{+}$. By \cite{Pal}, there exist $\alpha,\beta\in \mathbb{C}$ with $\beta^{d-1}=1$ and $\alpha^{d+1}=\beta$ such that $\beta p_c(y)=\alpha p_{c'}(\alpha y)$. Comparing coefficients, we get $c'=\alpha^d c$. Since $c,c'>0$, the quotient $c'/c=\alpha^d$ is a positive real number. Further $\alpha^d$ is a root of unity and only positive real root of unity is $1$. Therefore, $c'=c$.
It follows that if $c\neq c'$, then $I_{H_c}^{+}$ and $I_{H_{c'}}^{+}$ are not biholomorphic. The family $\{H_c:c\in\mathbb{R}_{>0}\}$ contains continuum many H\'enon maps of degree $d$ with mutually non-biholomorphic escaping sets. Moreover by \cite[Section~3.3]{BV}, it follows that if $d$ is prime, then $\mathrm{Aut(I_{H_c}^+)}\cong \mathrm{Aut}_1(I_{H_c}^+)\times \mathbb{Z}\cong \mathbb{C}\times\mathbb{Z}$.

\subsection*{Proof of Corollary \ref{global aut group}}
   Let $f\in\mathrm{Aut}_1(I_H^+)\cap\mathrm{Aut}(\mathbb{C}^2)$. It follows from \eqref{element in Aut escaping} that, $f=f_p^r\circ f_{\gamma'}$ for some $r\in\{1,2,\dots,i_p\}$ and $\gamma'\in \mathbb{C}$. Therefore $f_{\gamma'}\in \mathrm{Aut}_1(I_H^+)\cap\mathrm{Aut}(\mathbb{C}^2)$. It follows from \cite[Theorem 4.1]{BV}, that $f_{\gamma'}$ is identity map. Thus 
    \[
    \mathrm{Aut}_1(I_H^+)\cap\mathrm{Aut}(\mathbb{C}^2)=\{f_p,f_p^2,\dots,f_p^{i_p}\}\cong \mathbb{Z}_{i_p}.
    \]

\subsection*{Proof of Corollary \ref{level sets of green}}
It follows from \cite{BPV1}, that the restrictions of $\pi_H$, $\pi_H:\mathbb{C}\times \{t\in \mathbb{C}:1<|t|<e^c\} \to \Omega_c'$ is a covering map. Moreover, if $f \in \mathrm{Aut}_1(\Omega_c')$, then $f$ lifts to an automorphism of the corresponding covering space of the form $(s,t)\mapsto (\beta s+\gamma,\alpha t)$, where $\alpha^{d^2-1}=1$, $\beta=\alpha^{d+1}$, and $\gamma\in \mathbb{C}$. Furthermore, it follows from \cite[Equation (6.13)]{BPV1} that, the coefficients of $Q_H$ satisfy the same relations as those in \cite[Equation (4.13)]{Pal}. Consequently, the coefficients of the polynomial $p$ satisfy similar relations as in \eqref{alpha coeff poly}.  Since every automorphism of the escaping preserves the level sets of the Green's function, the automorphism $f_p$ constructed in Step~2 belongs to $\mathrm{Aut}_1(\Omega_c')\cap \mathrm{Aut}_1(\Omega^c)$. Therefore, using the same argument as above, we obtain $\mathrm{Aut}_1(\Omega'_c)\cong \mathbb{Z}_{i_p}\ltimes\mathbb{C}.$

\medskip
\no
It follows from the proof of corollary 1.3 in \cite{BV} that $\mathrm{Aut}(\Omega_c)\subset \mathrm{Aut}_1(I_H^+)\cap \mathrm{Aut}(\mathbb{C}^2)$. Conversely, if $f\in\mathrm{Aut}_1(I_H^+)\cap \mathrm{Aut}(\mathbb{C}^2)$, then, by \cite{MP}, $G_H^+\circ f=G_H^+$. Thus $f$ preserves the sub-level sets of $G_H^+$. In particular, $f(\Omega_c)=\Omega_c$. Thus $\mathrm{Aut}(\Omega_c)=\mathrm{Aut}_1(I_H^+)\cap \mathrm{Aut}(\mathbb{C}^2)=\{f_p,f_p^2,\dots, f_p^{i_p}\}$.

\subsection*{Proof of Corollary \ref{whole aut escaping}}
It follows from proof of Theorem 3.4 in \cite{BV} that $\mathrm{Aut}(I_H^+)\cong \mathrm{Aut}_1(I_H^+)\times \tilde{I}$, where $\mathbb{Z}\leq \tilde{I}\leq\mathbb{Z}^l$. Thus $\mathbb{Z}_{i_p}\ltimes\mathbb{C}\times\mathbb{Z}\leq \mathrm{Aut}(I_H^+)\leq \mathbb{Z}_{i_p}\ltimes\mathbb{C}\times\mathbb{Z}^l$.

\subsection*{Proof of Corollary \ref{global whole aut group}}
Let $f \in \mathrm{Aut}(I_H^+)\cap\mathrm{Aut}(\mathbb{C}^2)$ and suppose that the degree of $H$ is $d = p_1^{m_1}p_2^{m_2}\cdots p_l^{m_l}$, where $p_1,p_2,\dots,p_l$ are distinct primes and $m_1,m_2,\dots,m_l\in\mathbb{N}$. By \cite[Theorem 1.2]{MP}, there exist $s_1 > 0$, $s_2 \in \mathbb{Z}$, and a linear map $L$ such that
\[
f^{s_1}=L\circ H^{s_2} .
\]
The homomorphism $H_*:\mathbb{Z}[1/d]\to\mathbb{Z}[1/d]$ is given by $x\mapsto d x$, and the homomorphism $L_*:\mathbb{Z}[1/d]\to\mathbb{Z}[1/d]$ is identity map. Moreover, for some integers $n_1,n_2,\dots, n_l$, $f_*:\mathbb{Z}[1/d]\to\mathbb{Z}[1/d]$ is of the form $x\mapsto (-1)^{n_0}p_1^{n_1}p_2^{n_2}\cdots p_l^{n_l}x$. Consequently,
\[
\left((-1)^{n_0}p_1^{n_1}p_2^{n_2}\cdots p_l^{n_l}\right)^{s_1}=\left(p_1^{m_1}p_2^{m_2}\cdots p_l^{m_l}\right)^{s_2}.
\]
Comparing the exponents of the primes $p_i$, we obtain
$n_i s_1 \;=\; m_i s_2$ for $1\leq i \leq l.$ In particular, if $d=p_1p_2\cdots p_l$, then $n_i s_1=s_2$ for $1\leq i\leq l$. Therefore $n_1=n_2=\cdots=n_l$. Thus after composing $f$ with a suitable iterate of $H$, the resulting automorphism induces either $\mathrm{Id}$ or $-\mathrm{Id}$ on the fundamental group. As observed by Bousch \cite{Bousch}, $-\mathrm{Id}$ is not possible ( also see \cite[Theorem 3.4]{BV}). Hence, after this normalization, the resulting automorphism lifts to an automorphism of the covering space $\mathbb{C}\times\mathbb{C}\setminus\mathbb{\overline{D}}$.
Hence by \ref{global aut group}, $\mathrm{Aut}(I_H^+)\cap \mathrm{Aut}(\mathbb{C}^2)\cong \mathbb{Z}_{i_p}\times\mathbb{Z}$.

\subsection*{Proof of Theorem \ref{Biholo escaping}}
 By \cite[Equation 5.16]{Pal}, it follows that there exists $\alpha\in\mathbb{C}$ with $|\alpha|=1$ such that $a_{d-i}^H=\alpha a_{d-i}^F$, for $2\leq i\leq d$. Thus 
 \[
 gcd\{d^2-1,i: a^H_{d-i}\neq 0\}=gcd\{d^2-1,i: a^F_{d-i}\neq 0\}.
 \]
 Fix $f_\circ\in \mathrm{Biholo}_1(I_H^+\to I_F^+)$. Then for every $f\in \mathrm{Biholo}_1(I_H^+\to I_F^+)$ there exist $g\in \mathrm{Aut}_1(I_H^+)$ and $h\in\mathrm{Aut}_1(I_F^+)$ such that $f=f_\circ \circ g=h\circ f_\circ$. Therefore using \eqref{element in Aut escaping}, for every $f\in \mathrm{Biholo}_1(I_H^+\to I_F^+)$ there exist $r\in\{1,2,\dots,i_p\}$ and $\gamma\in \mathbb{C}$ such that 
\[
f=f_\circ \circ f_p^r\circ f_\gamma.
\]

\bibliographystyle{amsplain}

\end{document}